\documentclass[a4paper,11pt,reqno,twoside]{amsart}
\usepackage{ricotta_templier}
\begin{document}
\begin{merci}
Le premier auteur, financé par le projet  ANR <<Aspects Arithmétiques des Matrices Aléatoires et du Chaos Quantique>>, voudrait remercier chaleureusement H. Darmon pour lui avoir suggéré cette problématique et I. Fesenko pour son invitation et pour les excellentes conditions de travail offertes par l'Université de Nottingham. Le second auteur tient à remercier l'Université Bordeaux 1 pour son invitation au séminaire de théorie des nombres et pour son accueil chaleureux.
\end{merci}
\section{Introduction et présentation des résultats}\label{intro}%
\subsection{Problématique et résultats précédents}%
Soient $N$ un entier sans facteur carré et
\begin{equation*}
\mathcal{D}\coloneqq\left\{d\in\mathbb{Z},\;d<0,\;\mu^2(d)=1,\;d\equiv \nu^2 \mod (4N),\;(\nu,4N)=1\right\}
\end{equation*}
l'ensemble des discriminants fondamentaux impairs qui satisfont la \emph{condition de Heegner}. Pour $d$ dans $\mathcal{D}$, notons $\mathbb{H}_d$ le corps de classes de Hilbert du corps quadratique imaginaire $\Q(\sqrt{d})$. Soient $E$ une courbe elliptique rationnelle de conducteur $N$ et P$_d\in E(\mathbb{H}_d)$ l'un des points de Heegner de $E$ de discriminant $d$. Notons $\widehat{h}_{\mathbb{H}_d}($P$_d)$ la hauteur de Néron-Tate de ce point, sachant que la normalisation adoptée est donnée par $\widehat{h}_{\mathbb{H}_d}($P$_d)=[\mathbb{H}_d:\mathbb{Q}]\widehat{h}(\ptd)$ où $\widehat{h}:E(\overline{\mathbb{Q}})\rightarrow \mathbb{R}_+$ désigne la fonction hauteur canonique (voir \cite[Chapitre VIII Section 9]{MR1329092} par exemple).
Cette valeur, qui ne dépend pas du choix du point de Heegner de discriminant $d$ (voir \cite{GrZa}), est un invariant arithmétique du couple $(E,d)$ qui semble très irrégulier, comme le suggère la figure 1 de \cite{RiVi}.

Le but de cet article, est de poursuivre l'étude quantitative de $\widehat{h}_{\hd}(\ptd)$ \emph{en moyenne} sur les discriminants $d$ dans $\mathcal{D}$. De façon plus précise, il s'agit, comme dans \cite{RiVi}, d'étudier le comportement asymptotique de la somme
\begin{equation}
\label{eq_mohauteur}
\sum_{\substack{d\in\mathcal{D} \\
\abs{d}\leq Y}}\widehat{h}_{\mathbb{H}_d}(\text{P}_d)
\end{equation}
lorsque $Y>0$ tend vers l'infini.

Il est possible de montrer que l'ensemble $\mathcal{D}\subset \mathbb{Z}$ possède une densité 
\begin{equation*}
\lim_{Y\rightarrow+\infty}
\dfrac{%
\text{card}
\{
d\in \mathcal{D},\
|d|\le Y
\}
}%
{Y}
=2c_N
\end{equation*}
dont la valeur explicite est
\begin{equation*}
c_N\coloneqq\frac{3}{\pi^2N}\left(\prod_{\substack{p\in\mathcal{P} \\
p\mid 2N}} \left(1-\frac{1}{p^2}\right)^{-1}\right)\text{card}\left(\left\{\nu^2\in\mathbb{Z}\left/4N\mathbb{Z}\right., \;(\nu,4N)=1\right\}\right).
\end{equation*}
La plupart des notations de cet article sont compatibles avec celles de \cite{RiVi}.

L'approche analytique naturelle de cette étude repose sur la formule de Gross-Zagier (\cite[Théorème 6.1, Section I page 230]{GrZa})
\begin{equation}
\label{eq_gross_zagier}
L^\prime_d(E,1)=\frac{2\Omega}{u^2\sqrt{-d}}\widehat{h}_{\mathbb{H}_d}(\emph{P}_d)
\end{equation}
où $2u$ est le nombre de racines de l'unit\'e de $\Q(\sqrt{d})$, $\Omega$
est le volume complexe de $E$ c'est-à-dire le double de l'aire d'un parallélogramme fondamental de $E(\mathbb{C})$ (voir remarque \ref{rem_cP}). La définition de la série $L$ apparaissant dans le membre de gauche est donnée
dans la section~\ref{ana}.
Les auteurs de \cite{RiVi} ont déterminé le comportement asymptotique suivant (voir \cite[Théorème 4.1]{RiVi}).
\begin{theo}[G. Ricotta \& T. Vidick]
\label{traceana}
Si $E$ est une courbe elliptique rationnelle de conducteur $N$ sans facteur carré et $F$ est une fonction lisse à support compact dans $\mathbb{R}^\times_+$ alors
$$\sum_{d\in\mathcal{D}}L^\prime_d(E,1)F\left(\frac{\abs{d}}{Y}\right)=\alpha Y\log{Y}+\beta Y+\mathsf{Error}+O_\epsilon\left(N^{15/4+\epsilon}Y^{19/20+\epsilon}\right)$$
où
\begin{equation}
\label{conjerror}
\mathsf{Error}=O_\epsilon\left(NY\left(\log{(NY)}\right)^{1/2+\epsilon}\right)
\end{equation}
pour tout $\epsilon>0$ et o\`u
\begin{align*}
\alpha & \coloneqq c_NL(1)\int_0^{+\infty}F(t)\mathrm{d}t, \\
\beta & \coloneqq c_N\int_0^{+\infty}F(t)\left(L^\prime(1)+L(1)\left(\log{\left(\frac{Nt}{4\pi^2}\right)}-2\gamma\right)\right)\mathrm{d}t
\end{align*}
avec
\begin{equation*}
L(s)\coloneqq L(\mathrm{Sym}^2E,2s)\frac{\zeta^{(N)}(4s-2)}{\zeta^{(N)}(2s)}\prod_{(p,2N)=1}\left(1-\frac{1}{p^{4s-2}(p+1)}\right).
\end{equation*}
\end{theo}
\begin{remint}
Selon l'estimation \eqref{conjerror} du terme $\mathsf{Error}$, ce théorème est un développement asymptotique à un terme du premier moment de la valeur de la dérivée de la série $L_d(E,s)$ au point critique $s=1$, lorsque $F$ est d'intégrale non nulle. L'analyse analytique et numérique menée dans \cite{RiVi} avait suggéré que ce terme devrait être négligeable devant $Y$. L'objectif de ce travail est de confirmer cette observation c'est-à-dire de prouver qu'il existe des constantes réelles explicites $\delta_2$ et $\delta_1>0$ telles que
\begin{equation*}
\mathsf{Error}=O_\epsilon\left(N^{\delta_2+\epsilon}Y^{1-\delta_1+\epsilon}\right)
\end{equation*}
pour tout $\epsilon>0$ (voir \eqref{eq_goal}).
\end{remint}
\begin{remint}
La fonction $L(s)$ ci-dessus correspond à la fonction $\widetilde{\mathcal{L}}(s)$ de \cite{RiVi}. D'une part, son expression a été simplifiée: $L(s)$ apparaît clairement comme $L(\text{Sym}^2E,2s)$ multiplié par des facteurs correctifs, dont le produit converge absolument lorsque $s$ est proche de $1$. D'autre part, une petite erreur dans les facteurs Eulériens à la place $2$ a été corrigée. Ceci modifie très légèrement la valeur de $\beta$, mais ne modifie pas la valeur de $\alpha$. C'est uniquement ce facteur $\alpha$ qui est utilis\'e dans le paragraphe 5 de \cite{RiVi} lors de la confrontation des résultats numérique et théorique.
\end{remint}
%
\subsection{Présentation des nouveaux résultats}%
Le résultat principal de cet article est le suivant.
\begin{theoint}
\label{resultat_A}
Si $E$ est une courbe elliptique rationnelle de conducteur $N$ sans facteur carré et $F$ est une fonction lisse à support compact dans $\mathbb{R}^\times_+$ alors
\begin{equation*}
\sum_{d\in\mathcal{D}}L^\prime_d(E,1)F\left(\frac{\abs{d}}{Y}\right)=\alpha Y\log{Y}+\beta Y+O_\varepsilon\left(N^{15/4+\epsilon}Y^{19/20+\epsilon}\right)
\end{equation*}
pour tout $\varepsilon>0$.
\end{theoint}

\begin{remint} L'estim\'ee est tr\`es pr\'ecise en la variable $Y$. Si la courbe elliptique $E$ est fix\'ee, une petite puissance dans le terme d'erreur a été sauvée par rapport au terme principal. Bien qu'il soit possible d'am\'eliorer nettement les exposants $15/4$ et $19/20$, des arguments courts, qui donnent des exposants qui ne sont pas optimaux mais restent de taille mod\'er\'ee, ont été privil\'egi\'es. 

Une dépendance explicite en le conducteur $N$ a été conservé tout au long de ce travail. La façon dont varie $\widehat{h}_{\mathbb{H}_d}(\emph{P}_d)$ (ou plus simplement $\widehat{h}(E)$, la hauteur de Faltings de $E$) en fonction de $N$ est une question d\'elicate et largement ouverte.
\end{remint}

\begin{remint}
D'une part, $c_N\asymp 2^{-\omega(N)}$ et $L(1)\asymp L(\mathrm{Sym}^2 E,2)$ . D'autre part,
\begin{equation*}
 N^{-\epsilon} \ll_\epsilon L(\mathrm{Sym}^2 E,2) \ll \log N
\end{equation*}
d'après \cite{GoHoLi}. Ainsi, le terme d'erreur croît (à $Y$ fixé) bien plus vite avec $N$ que les termes principaux ($\alpha$ et $\beta$). Le terme principal domine lorsque $N \ll Y^{1/76}$.
\end{remint}
Le th\'eor\`eme~\ref{resultat_A} et la formule de Gross-Zagier \eqref{eq_gross_zagier} impliquent le corollaire suivant (voir \cite{RiVi}).
\begin{corint}
Si $E$ est une courbe elliptique rationnelle de conducteur $N$ sans facteur carré alors
$$\sum_{\substack{d\in\mathcal{D} \\
\abs{d}\leqslant Y}}\widehat{h}_{\mathbb{H}_d}(\emph{P}_d)=C_{\emph{P}}Y^{\frac{3}{2}}\log{Y}+C_{\emph{P}}^\prime Y^{\frac{3}{2}}+O_\epsilon\left(N^{15/4+\epsilon}Y^{29/20+\epsilon}\right)$$
pour tout $\varepsilon>0$, o\`u $C_{\emph{P}}$ est la constante d\'efinie par
\begin{equation*}
  C_{\emph{P}}\coloneqq\frac{\pi}{3}c_N\prod_{\substack{p\in\mathcal{P} \\
  (p,2N)=1}}\left(1-\frac{1}{p^2(p+1)}\right)^{-1}\frac{L(\emph{Sym}^2 E,2)}{\pi\Omega}
\end{equation*}
et
\begin{equation*}
C_{\emph{P}}^\prime\coloneqq C_{\emph{P}}\left(\log{\left(\frac{N}{4\pi^2}\right)}-\frac{2}{3}-2\gamma\right)+\frac{c_N}{3\Omega}L^\prime(1).
\end{equation*}
\end{corint}
\begin{remint}\label{rem_cP}
Une autre écriture possible pour la constante $C_{\emph{P}}$ est
\begin{equation*}
C_{\emph{P}}=\frac{\pi}{3}c_N\prod_{\substack{p\in\mathcal{P} \\
  (p,2N)=1}}\left(1-\frac{1}{p^2(p+1)}\right)^{-1}\frac{\text{deg}(\Phi)}{N}
\end{equation*}
où $\Phi:X_0(N)\rightarrow E$ est la paramétrisation modulaire de $E$. Cela résulte essentiellement de \cite{MR0434962} et de l'égalité $\Omega\times\text{deg}(\Phi)=4\pi^2(f,f)$ où $(f,f)$ désigne le produit scalaire de Peterson de $f$ (voir \cite[Pages 230, 308 et 310]{GrZa}). Ainsi, l'intuition géométrique suggérant que, à conducteur $N$ fixé, la hauteur de Néron-Tate d'un point de Heegner est gouvernée par le degré de la paramétrisation modulaire de $E$, est vraie en moyenne sur les discriminants (voir \cite{RiVi}).
\end{remint}
\subsection{Idée de la preuve du Théorème \ref{resultat_A}}%
Les premi\`eres \'etapes sont classiques et ont \'et\'e effectu\'ees dans \cite{RiVi}. On exprime $L'_d(E,1)$ comme un polyn\^ome de Dirichlet tronqu\'e \`a l'aide de l'\'equation fonctionnelle approchée. Un terme <<diagonal>> provenant des entiers $n$ qui sont des carr\'es fournit le terme principal $\alpha Y\log Y +\beta Y$, avec $\alpha$ et $\beta$ comme données pr\'ecedemment (not\'e $\mathsf{TP}_1$ dans \cite{RiVi}). La contribution <<hors-diagonale>> $r'_d(n)$, voir la d\'efinition~\eqref{eq_r_d}, correspond aux entiers $v$ qui sont diff\'erents de $0$. C'est d'elle que provient le terme $\mathsf{Error}$. Pour d\'emontrer le Th\'eor\`eme~\ref{resultat_A}, on améliore l'estim\'ee \eqref{conjerror}. L'estim\'ee~\eqref{eq_goal} montre  que $\mathsf{Error}$ est n\'egligeable devant $\beta Y$ et que l'on peut l'incorporer dans le terme d'erreur.

Expliquons la différence principale avec l'article ant\'erieur. Dans \cite{RiVi}, il est fait appel \`a \cite[(4.3)]{RiVi}:
\begin{equation}\label{eq_rivi}
  \sum_{d\in \mathcal{D}} r'_d(n) \leq 4 \sqrt{n}.
\end{equation}
Cette majoration combine astucieusement les aspects <<clairsemés>> de l'ensemble $\mathcal{D}$ et de l'équation $4n=u^2+\abs{d}v^2$. Elle met en \'evidence des annulations non triviales dans la somme en $d$: il en d\'ecoule l'estim\'ee \eqref{conjerror}.

Dans cet article, les annulations de la somme en $d$ sont davantage exploitées en prenant en compte les oscillations des coefficients de Fourier $a_n$ dans les progressions arithmétiques (qui ne peuvent pas être d\'etect\'ees avec \eqref{eq_rivi}). Chaque étape est donc délicate dans la mesure o\`u toute majoration trop directe briserait ces oscillations. C'est pour cette raison que la premi\`ere \'etape de la d\'emonstration consiste \`a <<d\'eployer>> l'\'equation $4n=u^2+\abs{d}v^2$.  Intuitivement, lorsque $d$ parcourt $\mathcal{D}$, l'entier $n$ parcourt des progressions arithm\'etiques. Au coeur de la démonstration est l'inégalité \eqref{eq_estimcle}.

L'estimée \eqref{eq_estimcle} est très g\'en\'erale et peut être améliorée dans de nombreux cas \cite{MR1111014,MR1111010,MR1220457,MR1476732,MR1262431,MR0126423,MR0140491,MR1039952}. Elle est suffisante dans le cas pr\'esent car le module des progressions arithmétiques est de taille mod\'er\'ee. La d\'emonstration est de ce fait robuste et il est probable qu'elle s'adapte \`a d'autres problèmes.
\subsection{Organisation de l'article}%
La Section \ref{ana} d\'ecrit les propriétés analytiques de la série $L$ en \eqref{eq_gross_zagier}, ce qui permet de décrire précisément le terme $\mathsf{Error}$. La d\'emonstration du Théorème \ref{resultat_A} est détaillée dans la Section~\ref{preuve}.
\begin{notations}
Le paramètre principal de cet article est un nombre réel strictement positif $\;Y$ qui tend vers l'infini. Ainsi, si $f$ et $g$ sont des fonctions de la variable réelle à valeurs complexes alors les notations $f(Y)\ll_{B}g(Y)$ ou $f(Y)=O_B(g(Y))$ signifient que $\abs{f(Y)}$ est inférieur ou égal à une <<constante>> ne dépendant que de $B$ multipliée par $g(Y)$, au moins pour $Y>0$ assez grand. La lettre $\epsilon$ désigne un réel strictement positif arbitrairement petit dont la valeur peut changer d'une ligne à l'autre. La lettre $\mu$ désigne la fonction de Möbius et la lettre $\omega$ désigne la fonction nombre de diviseurs premiers. Enfin, pour un produit Eulerien $L(s)=\prod_{p\in\mathcal{P}}L_p(s)$, $\mathcal{P}$ désignant l'ensemble des nombres premiers, et pour un entier $k\geq 1$, posons
\begin{equation*}
L_{(k)}(s)\coloneqq\prod_{\substack{p\in\mathcal{P} \\
p\mid k}}L_p(s) \quad\text{ et } \quad L^{(k)}(s)\coloneqq\prod_{\substack{p\in\mathcal{P} \\
(p,k)=1}}L_p(s).
\end{equation*}
\end{notations}
%
\section{Prérequis analytiques}\label{ana}%
La fonction $L$ de $E$ sur $\mathbb{Q}$ notée $L(E\vert\mathbb{Q},s)=\sum_{n\geqslant 1}a_n n^{-s}$ est définie a priori sur $\Re{(s)}>\frac{3}{2}$. Les travaux de A. Wiles et R. Taylor (\cite{Wi}, \cite{TaWi}) assurent qu'il existe une forme primitive cuspidale $f$ de niveau $N$, de poids $2$ et de caractère trivial telle que
\begin{equation*}
L(E\vert\mathbb{Q},s)=L(f,s).
\end{equation*}
Souvenons-nous que les coefficients de Fourier de $f$ satisfont
\begin{equation}\label{eq_deligne}
a_n\ll_\epsilon n^{1/2+\varepsilon}
\end{equation}
pour tout $\varepsilon>0$. La série de Dirichlet apparaissant dans le membre de gauche de la formule de Gross-Zagier \eqref{eq_gross_zagier} est définie sur $\Re{(s)}>3/2$ par
\begin{equation*}
L_d(E,s)\coloneqq\left(\sum_{\substack{m\geqslant 1 \\
(m,N)=1}}\frac{\chi_d(m)}{m^{2s-1}}\right)\times\left(\sum_{n\geqslant 1}\frac{a_n r_d(n)}{n^s}\right)
\end{equation*}
où $\chi_d$ est le caractère de Kronecker du corps $\Q(\sqrt{d})$ et $r_d(n)$ désigne le nombre d'idéaux principaux de $\Q(\sqrt{d})$ de norme $n$. Par la méthode de Rankin-Selberg, $L_d(E,s)$ admet un prolongement holomorphe à $\mathbb{C}$ et satisfait l'\'equation fonctionnelle
$$\forall s\in\mathbb{C}, \quad \Lambda_d(E,s)=-\chi_d(N)\Lambda_d(E,2-s)$$
où $\Lambda_d(E,s)\coloneqq(N\abs{d})^s\left((2\pi)^{-s}\Gamma(s)\right)^2L_d(E,s)$ est la série complétée (voir la section IV de \cite{GrZa} pour la mise en oeuvre de la méthode). Remarquons que comme $d$ est un carré modulo $N$, le signe de l'équation fonctionnelle est $-1$, de sorte que $L_d(E,1)=0$.

Comme expliqué dans l'introduction, il s'agit d'estimer minutieusement le terme $\mathsf{Error}$ qui est donné dans les équations (4.4) et (4.5) de \cite{RiVi}, et qui apparaît lorsque l'on applique \`a $\Lambda_d(E,s)$ la méthode de l'équation fonctionnelle approchée:
\begin{equation}\label{eq_error1}
\mathsf{Error}\coloneqq2\sum_{d\in\mathcal{D}}\sum_{(m,N)=1}\sum_{n\geqslant 1}\frac{a_{n}\chi_d(m)r_d^\prime(n)}{mn}V\left(\frac{4\pi^2 n m^2}{N\abs{d}}\right)F\left(\frac{\abs{d}}{Y}\right).
\end{equation}
Les notations sont les suivantes:
\begin{equation}\label{eq_r_d}
r_d^\prime(n)\coloneqq\text{card}\left(\left\{(u,v)\in\left(\N\times\Z^*\right),\;u^2+\abs{d}v^2=4n\right\}\right)
\end{equation}
pour tout entier naturel $n\geq 1$ et tout élément $d$ de $\mathcal{D}$ (et plus généralement $\mathcal{D}^\prime$ défini en \eqref{eq_Dprime}); $V:\R_+^\ast\rightarrow\R$ est une fonction de coupure définie à la page 14 de \cite{RiVi} vérifiant
\begin{equation}\label{eq_V}
x^jV^{(j)}(x)\ll_j x^{-1/4}\exp{(-2\sqrt{x})}
\end{equation}
pour tout entier naturel $j$.
\section{Estimée du terme $\mathsf{Error}$}%
L'objectif de cette section est de démontrer que
\begin{equation}\label{eq_goal}
  \mathsf{Error}=O_\varepsilon\left(N^{13/12+\epsilon}Y^{5/6+\epsilon}\right)
\end{equation}
pour tout $\epsilon>0$.
\label{preuve}%
\subsection{Estimées triviales}%
On commence par déployer la définition~\eqref{eq_r_d} de sorte que:
\begin{equation*}
\mathsf{Error}
=
2
\sum_{(m,N)=1}
\sum_{(d,n,u,v)}
\frac{a_{n}\chi_d(m)}
{mn}
V\left(\frac{4\pi^2 n m^2}{N\abs{d}}\right)
F\left(\frac{\abs{d}}{Y}\right)
\end{equation*}
où la deuxième somme porte sur l'ensemble des quadruplets $(d,n,u,v)\in \mathcal{D}\times \mathbb{N}^*\times \mathbb{N} \times \mathbb{Z}^*$ qui vérifient $4n=u^2+\abs{d}v^2$.

Le terme $F(\abs{d}\left/Y\right.)$ restreint la sommation à $Y \ll\abs{d}\ll Y$ car $F$ est à support compact dans $\R_+^\times$. La décroissance exponentielle de $V$ donnée en \eqref{eq_V} assure que la contribution des quadruplets pour lesquels \[nm^2\gg_{\epsilon} N^{1+\epsilon}Y^\epsilon\abs{d}\] est négligeable car elle est bornée par $O_{B,\epsilon}((NY)^{-B})$ pour tout nombre réel $B>0$. C'est en particulier le cas lorsque $u>U\left/m\right.$ ou bien
$\abs{v}>V\left/m\right.$ ou encore $n>N_0\left/m^2\right.$ où
\begin{equation*}
U\coloneqq N^{1/2+\epsilon/2}Y^{1/2+\epsilon/2},\quad
V\coloneqq N^{1/2+\epsilon/2}Y^{\epsilon/2},\quad
N_0\coloneqq (NY)^{1+\epsilon}.
\end{equation*}
Dans toute la suite, il est donc légitime de se restreindre aux intervalles de sommation
\begin{equation}\label{eq_restriction}
0 \leq u\leq \frac{U}{m}
\quad
\text{  et  }
\quad
1\leq\abs{v}\leq \frac{V}{m}
\quad
\text{  et   }
\quad
1 \leq m \leq V
\quad
\text{  et   }
\quad
1 \leq n\leq \frac{N_0}{m^2}.
\end{equation}
M\^eme lorsqu'elles ne sont pas explicit\'ees, ces in\'egalit\'es sont sous-entendues dans toute la suite.

\begin{remark}\label{rem_trivial}
  Comme $n$ est uniquement déterminé par $u,v$ et $d$, le reste de la contribution peut se majorer, en valeur absolue et à l'aide de~\eqref{eq_deligne}, par
\begin{equation}\label{rem_triviale}
\begin{aligned}
  \mathsf{Error}
  & \ll
  \sum_{m\geq 1}
  \sum_{(u,v,d)}
   \frac{1}{m(u^2+\abs{d}v^2)^{1/2}}\left(\frac{N\abs{d}}{(u^2+\abs{d}v^2)m^2}\right)^{1/4} \\
  & \ll N^{1/4}\sum_{m\geq 1} \frac{1}{m^{3/2}}\left\{\sum_{\substack{\abs{d}\ll Y \\
  1\leq\abs{v}\leq V/m}}\frac{1}{\abs{d}^{1/2}\abs{v}^{3/2}}+\sum_{\substack{\abs{d}\ll Y \\
  1\leq u\leq U/m \\
  1\leq\abs{v}\leq V/m}}\frac{1}{u\abs{v}^{1/2}}\right\} \\
  & \ll_\epsilon N^{1/2+\epsilon}Y^{1+\epsilon}
 \end{aligned}
 \end{equation}
pour tout $\epsilon>0$. Cette estimée <<triviale>> n'est pas suffisamment précise (en fait moins bonne que \eqref{conjerror}). Pour obtenir~\eqref{eq_goal}, il faut tenir compte des oscillations des coefficients $a_n$.
\end{remark}
\subsection{La condition $d$ sans facteur carré}%
Pour effectuer la somme sur $d\in \mathcal{D}$, il est commode d'introduire l'ensemble
\begin{equation}\label{eq_Dprime}
\mathcal{D}^\prime\coloneqq
\left\{
d\in\mathbb{Z},\;d<0,\;d\equiv\nu^2\mod (4N),\;(\nu,4N)=1
\right\}.
\end{equation}
qui est une union de $O(N)$ progressions arithmétiques de raison $4N$; et de détecter séparément la condition $d$ sans facteur carré à l'aide de l'identité $\sum_{a^2\mid d}\mu(a)=\mu^2(d)$. La définition~\eqref{eq_r_d} est étendue à $d\in\mathcal{D}'$ et la somme est coupée selon la taille des diviseurs $a$:
\begin{equation*}
\mathsf{Error} \leq \sum_{1 \leq a \leq A}
\abs{\mathsf{Error}(a)} +\sum_{a>A} \abs{\mathsf{Error}(a)} = E_1+E_2.
\end{equation*}
Ici, $A$ est une petite puissance de $Y$ qui sera <<optimisée>> dans la section \ref{sec_resultat_A} et
\begin{equation*}
\mathsf{Error}(a)
\coloneqq
2\sum_{
\substack{%
d\in\mathcal{D}^\prime \\
a^2|d
}}
\sum_{(m,N)=1}
\sum_{n\geqslant 1}
\frac{a_{n}\chi_{d}(m)r_{d}^\prime(n)}{mn}
V\left(\frac{4\pi^2 n m^2}{N\abs{d}}\right)
F\left(\frac{\abs{d}}{Y}\right).
\end{equation*}
Le terme $\mathsf{E_1}$ est estimé dans la section suivante alors que le terme $\mathsf{E_2}$ est estimé dans la section~\ref{large}.
\subsection{Contribution des petits diviseurs}\label{small}%
On suppose que $(a,4N)=1$ car sinon $\mathsf{Error}(a)=0$.
\subsubsection{Travail préparatoire}
La loi de réciprocité quadratique implique que
\begin{equation*}
\chi_d(m)=\begin{cases}
\legendre{d}{m_1} & \text{si $(m_2,d)=1$ et $m_1$ impair,} \\
\chi_8(d) \legendre{d}{m_1} & \text{si $(m_2,d)=1$ et $m_1$ est pair} \\
0 & \text{si $(m_2,d)\not=1$}
\end{cases}
\end{equation*}
pour tout $d$ dans $\mathcal{D}^\prime$ et tout entier naturel $m=m_1m_2^2\geq 1$ où $m_1$ est sans facteur carré. En particulier, la fonction
\begin{equation*}
d\mapsto \chi_d(m)
\end{equation*}
est $4m$-périodique.

Supposons, afin d'alléger les notations, que $N$ est impair. Le cas $N$ pair est similaire. Pour tout entier $d<0$, l'appartenance $d\in\mathcal{D}^\prime$ peut être détectée par l'identité:
\begin{equation*}
\frac{1}{2^{\omega(N)+1}}\sum_{\chi\!\!\!\!\!\!\!\pmod{4}}\sum_{\substack{\psi\!\!\!\!\!\!\!\pmod{N} \\
\psi^2=1}}
\chi(d)\psi(d)
=
\delta_{d\in\mathcal{D}^\prime}
\end{equation*}
où $\chi$ décrit les deux caractères de Dirichlet de module $4$ et $\psi$ décrit les caractères de Dirichlet quadratiques de module $N$.

Pour toute la suite, posons
\begin{equation*}
  \phi(d)\coloneqq \chi_d(m)\delta_{d\in\mathcal{D}'}.
\end{equation*}
C'est une fonction $4mN$-périodique.

En tenant compte de \eqref{eq_r_d}, il est possible d'obtenir l'estimation
\begin{equation}\label{eq_error2}
\mathsf{Error}(a)
\leq
\sum_{(m,N)=1}
\frac{1}{m}
\sum_{(u,v)\in \N\times \Z^\ast}
\left
\vert
\sum_{(d,n)}
\phi(d)
\frac{a_n}{n}
V\left(\frac{4\pi^2 n m^2}{N\abs{d}}\right)
F\left(\frac{\abs{d}}{Y}\right)
\right
\vert
\end{equation}
où la sommation porte sur les couples $(d,n)\in \mathbb{Z}_-\times \mathbb{N}^\ast$ satisfaisant $4n=u^2+\abs{d}v^2$, $a^2\mid d$ et $1 \leq n\leq N_0/m^2$.

Limitons-nous à estimer la contribution des $v$ pairs sachant que la contribution des $v$ impairs se traite de façon tout à fait similaire. Dans ce cas, $u$ est également pair. Il s'agit alors d'estimer
\begin{equation}\label{eq_error3}
\sum_{\substack{(m,N)=1 \\
1\leq m\leq V}}
\frac{1}{m}\sum_{\substack{%
0\leq u\leq U\left/(2m)\right. \\
1\leq\abs{v}\leq V\left/(2m)\right.}}
\left\vert\mathsf{E}(u^2;a^2,v^2)\right\vert,
\end{equation}
après avoir effectué le changement de variables $(u,v)\mapsto(u/2,v/2)$, et où
\begin{equation}\label{eq_prog}
\mathsf{E}
(u^2;a^2,v^2)
\coloneqq
\sum_{\substack{%
n\equiv u^2\pmod{a^2v^2} \\
u^2+a^2v^2\leq n\leq N_0\left/m^2\right.}}
\frac{a_{n}}{n}
\phi
\left(\frac{u^2-n}{a^2v^2}\right)
V\left(\frac{4\pi^2 n m^2}{N\left(n-u^2\right)\left/v^2\right.}\right)
F\left(\frac{(n-u^2)\left/v^2\right.}{Y}\right)
\end{equation}
\subsubsection{Coefficients de Fourier dans des progressions arithmétiques de petit module}%
Selon \eqref{eq_prog}, il s'agit donc d'estimer
\begin{equation*}
\mathsf{E}(u^2;a^2,v^2)=
\sum_{n_1 \leq n \leq n_2}
a_{n}
\eta(n)
G(n),
\end{equation*}
o\`u $1\leq  n_1\coloneqq u^2+a^2v^2\ll n_2\coloneqq N_0/m^2$, $\eta$ est la fonction $4Nma^2v^2$-périodique et de module inférieur à $1$ définie par
\begin{equation*}
  \eta(n)
\coloneqq
\begin{cases}
\phi\left(\dfrac{u^2-n}{a^2v^2}\right) & \text{si $n\equiv u^2 \pmod{a^2 v^2}$,}\\
0 & \text{sinon,}
\end{cases}
\end{equation*}
et $G$ est la fonction définie par
\begin{equation*}
G(x)\coloneqq \frac{1}{x}V\left(\frac{4\pi^2 x m^2}{N\left(x-u^2\right)\left/v^2\right.}\right)F\left(\frac{\left(x-u^2\right)\left/v^2\right.}{Y}\right)
\end{equation*}
pour tout nombre réel $x>0$. Selon \eqref{eq_V},
\begin{equation}\label{eq_G}
\begin{aligned}
G(x)
& \ll
\left(\frac{NY}{m^2}\right)^{1/4}x^{-5/4}, \\
xG^\prime(x)
& \ll
\left(\frac{NY}{m^2}\right)^{1/4}x^{-5/4}
\left(
1
+
\frac{u^2}{Y^{1/4}\abs{v}^{1/2}(x-u^2)^{3/4}}
\right).
\end{aligned}
\end{equation}

D'après \cite[Lemma I]{MR1081731}, on a
\begin{equation}\label{eq_estimcle}
\sum_{1\leq n\leq x}
a_{n}
\eta(n)
\ll
(Nma^2v^2)^{1/2}
x \log x.
\end{equation}
En sommant par parties,
\begin{equation*}
\mathsf{E}(u^2;a^2,v^2)=
\left[
\left\{\sum_{1\leq n\leq x}a_{n}\eta(n)\right\}
G(x)
\right]_{n_1}^{n_2}
-
\int_{n_1}^{n_2}
\left\{
\sum_{1\leq n\leq x}
a_{n}\eta(n)
\right\}
G^\prime(x)\dd x
\end{equation*}
ce qui implique que
\begin{equation}\label{eq_inter}
\mathsf{E}(u^2;a^2,v^2)\ll_\epsilon
(NY)^\epsilon
\left(Nma^2v^2\right)^{1/2}
\left(\frac{NY}{m^2}\right)^{1/4}
\frac{N^{1/4}}
{\left(u^2+a^2v^2\right)^{1/4}}
\end{equation}
pour tout $\epsilon>0$ grâce à \eqref{eq_estimcle} et \eqref{eq_G}.
\subsubsection{Estimée finale de la contribution des petits diviseurs}%
Une estimation triviale des sommes en $u$, $v$, et $m$ combinée avec \eqref{eq_error3} et \eqref{eq_inter} assure que
\begin{equation}
\begin{split}
 \mathsf{Error}(a)
 &\ll_\epsilon
(NY)^\epsilon NY^{1/4}\sum_{m\leq V}\frac{1}{m}\left\{\sum_{1\leq\abs{v}\leq V\left/(2m)\right.}\sqrt{a}\sqrt{\abs{v}}+\sum_{\substack{1\leq u\leq U\left/(2m)\right. \\
1\leq\abs{v}\leq V\left/(2m)\right.}}\frac{a\abs{v}}{\sqrt{u}}\right\}
\\
& \ll_\epsilon (NY)^\epsilon NY^{1/4}\sum_{m\leq V}\frac{1}{m^2}\left\{\sqrt{a}V^{3/2}+a\sqrt{U}V^2\right\} \\
& \ll_\epsilon (NY)^\epsilon N^{9/4}Y^{1/2}a
\end{split}
\end{equation}
pour tout $\epsilon>0$ d'où
\begin{equation}\label{eq_small}
\mathsf{E_1} \ll_\epsilon (NY)^\epsilon N^{9/4}Y^{1/2}A^2.
\end{equation}
pour tout $\epsilon>0$. Cela conclut l'estimation de $\mathsf{E_1}$.
\subsection{Contribution des grands diviseurs}\label{large}%
En procédant comme dans la remarque~\ref{rem_trivial}, mais en tenant compte de la condition supplémentaire $a^2\mid d$, la majoration
\begin{equation*}
  \mathsf{Error}(a) \ll_\epsilon N^{1/2+\epsilon} Y^{1+\epsilon}/a^2
\end{equation*}
est valide pour tout $\epsilon>0$ et implique, en sommant sur les $a>A$, que
\begin{equation}\label{eq_large}
  \mathsf{E_2}\ll_\epsilon(NY)^\epsilon\frac{N^{1/2}Y}{A}
\end{equation}
pour tout $\epsilon>0$.
\subsection{Choix du paramètre $A$}\label{sec_resultat_A}%
Selon \eqref{eq_small} et \eqref{eq_large},
\begin{equation*}
\mathsf{Error}\ll_\epsilon (NY)^\epsilon\left(N^{9/4}Y^{1/2}A^2+\frac{N^{1/2}Y}{A}\right)
\end{equation*}
pour tout $\epsilon>0$ et le choix optimal est $A\coloneqq Y^{1/6}N^{-7/12}$.

\bibliographystyle{amsplain}
\bibliography{ricotta_templier}

\providecommand{\bysame}{\leavevmode\hbox to3em{\hrulefill}\thinspace}
\providecommand{\MR}{\relax\ifhmode\unskip\space\fi MR }
\providecommand{\MRhref}[2]{%
  \href{http://www.ams.org/mathscinet-getitem?mr=#1}{#2}
}
\providecommand{\href}[2]{#2}
\begin{thebibliography}{10}

\bibitem{MR0126423}
K.~Chandrasekharan and Raghavan Narasimhan, \emph{Hecke's functional equation
  and the average order of arithmetical functions}, Acta Arith. \textbf{6}
  (1960/1961), 487--503. \MR{MR0126423 (23 \#A3719)}

\bibitem{MR0140491}
\bysame, \emph{Functional equations with multiple gamma factors and the average
  order of arithmetical functions}, Ann. of Math. (2) \textbf{76} (1962),
  93--136. \MR{MR0140491 (25 \#3911)}

\bibitem{MR1111010}
W.~Duke and H.~Iwaniec, \emph{Estimates for coefficients of {$L$}-functions.
  {I}}, Automorphic forms and analytic number theory (Montreal, PQ, 1989),
  Univ. Montr\'eal, Montreal, QC, 1990, pp.~43--47. \MR{MR1111010 (92f:11068)}

\bibitem{MR1220457}
\bysame, \emph{Estimates for coefficients of {$L$}-functions. {II}},
  Proceedings of the Amalfi Conference on Analytic Number Theory (Maiori, 1989)
  (Salerno), Univ. Salerno, 1992, pp.~71--82. \MR{MR1220457 (94f:11041)}

\bibitem{MR1476732}
\bysame, \emph{Estimates for coefficients of {$L$}-functions. {III}},
  S\'eminaire de Th\'eorie des Nombres, Paris, 1989--90, Progr. Math., vol.
  102, Birkh\"auser Boston, Boston, MA, 1992, pp.~113--120. \MR{MR1476732
  (98h:11056)}

\bibitem{MR1262431}
\bysame, \emph{Estimates for coefficients of {$L$}-functions. {IV}}, Amer. J.
  Math. \textbf{116} (1994), no.~1, 207--217. \MR{MR1262431 (95k:11073)}

\bibitem{GoHoLi}
D.~Goldfeld, J.~Hoffstein, and Lieman D., \emph{An effective zero free region},
  Ann. of Math. (2) \textbf{140} (1994), no.~2.

\bibitem{GrZa}
Benedict~H. Gross and Don~B. Zagier, \emph{Heegner points and derivatives of
  {$L$}-series}, Invent. Math. \textbf{84} (1986), no.~2, 225--320.
  \MR{MR833192 (87j:11057)}

\bibitem{MR1039952}
James~Lee Hafner and Aleksandar Ivi{\'c}, \emph{On sums of {F}ourier
  coefficients of cusp forms}, Enseign. Math. (2) \textbf{35} (1989), no.~3-4,
  375--382. \MR{MR1039952 (90k:11055)}

\bibitem{MR1081731}
Henryk Iwaniec, \emph{On the order of vanishing of modular {$L$}-functions at
  the critical point}, S\'em. Th\'eor. Nombres Bordeaux (2) \textbf{2} (1990),
  no.~2, 365--376. \MR{MR1081731 (92h:11040)}

\bibitem{MR1111014}
R.~A. Rankin, \emph{Sums of cusp form coefficients}, Automorphic forms and
  analytic number theory (Montreal, PQ, 1989), Univ. Montr\'eal, Montreal, QC,
  1990, pp.~115--121. \MR{MR1111014 (92g:11042)}

\bibitem{RiVi}
Guillaume Ricotta and Thomas Vidick, \emph{Hauteur asymptotique des points de
  heegner}, à paraitre au Canadian Journal of Mathematics et disponible sur le
  site web de prépublication
  \url{http://journals.cms.math.ca/prepub/precjm.html} et à
  \url{http://www.math.u-bordeaux.fr/~ricotta/heegner.htm}.

\bibitem{MR0434962}
Goro Shimura, \emph{The special values of the zeta functions associated with
  cusp forms}, Comm. Pure Appl. Math. \textbf{29} (1976), no.~6, 783--804.
  \MR{MR0434962 (55 \#7925)}

\bibitem{MR1329092}
Joseph~H. Silverman, \emph{The arithmetic of elliptic curves}, Graduate Texts
  in Mathematics, vol. 106, Springer-Verlag, New York, 1992, Corrected reprint
  of the 1986 original. \MR{MR1329092 (95m:11054)}

\bibitem{TaWi}
Richard Taylor and Andrew Wiles, \emph{Ring-theoretic properties of certain
  {H}ecke algebras}, Ann. of Math. (2) \textbf{141} (1995), no.~3, 553--572.
  \MR{MR1333036 (96d:11072)}

\bibitem{Wi}
Andrew Wiles, \emph{Modular elliptic curves and {F}ermat's last theorem}, Ann.
  of Math. (2) \textbf{141} (1995), no.~3, 443--551. \MR{MR1333035 (96d:11071)}

\end{thebibliography}
\strut\newline
\noindent{\textit{G. Ricotta}\newline}
Université Bordeaux 1, Institut de Mathématiques de Bordeaux, Laboratoire A2X, 351, cours de la libération, 33405 Talence Cedex, France; \url{Guillaume.Ricotta@math.u-bordeaux1.fr}\newline
\noindent{\textit{N. Templier}\newline}
Université Montpellier 2, Institut de Mathématiques et de Modélisation de Montpellier, Case courrier 051, Place Eugène Bataillon, 34095 Montpellier Cedex, France; \url{templier@math.univ-montp2.fr}
\end{document}